\documentclass[reqno]{amsart}
\usepackage{amsfonts}
\usepackage{amsmath}
\usepackage{amsthm}
\usepackage[utf8]{inputenc}
\usepackage[T1]{fontenc}
\usepackage{graphicx} 
\usepackage{multicol}
\usepackage{caption}
\usepackage{color}
\usepackage{subcaption}
\usepackage{hyperref}
\usepackage{mathtools}
\usepackage{thmtools}
\usepackage{thm-restate}
\usepackage{cleveref}
\newtagform{brackets}{[}{]}
\usetagform{brackets}
\usepackage{soul}
\usepackage{float}
\allowdisplaybreaks
\usepackage{etoolbox}
\newcounter{temp}
\theoremstyle{plain}

\theoremstyle{definition}
\newtheorem*{Def}{Definition}

\newcommand{\figref}[1]{Figure~{\bf\ref{#1}}}

\newcommand{\secref}[1]{Section~{\bf\ref{#1}}}

\newcommand{\twoeqandorref}[2]{Equations~{\bf\ref{#1}}~and/or~{\bf\ref{#2}}}
\newcommand{\twoeqref}[2]{Equations~{\bf\ref{#1}}~and~{\bf\ref{#2}}}

\newcommand{\exref}[1]{Example~{\bf\ref{#1}}}

\everymath{\displaystyle}

\begin{document}
\title{STUDENT INQUIRY AND THE RASCAL TRIANGLE}
\author{Philip~K. Hotchkiss}
\address{Department of Mathematics\\Westfield State University\\Westfield, MA 01085}
\email{photchkiss@westfield.ma.edu}
\thanks{The author's work as part of the \textit{Discovering the Art of Mathematics Project} (\url{http://artofmathematics.org}), which provides the basis for the curriculum in this course was supported by NSF grants DUE-0836943, DUE-1225915, and a generous gift from Mr. Harry Lucas.}
\maketitle

\begin{abstract}
Those of us who teach Mathematics for Liberal Arts (MLA) courses often underestimate the mathematical abilities of the students enrolled in our courses.  Despite the fact that many of these students suffer from math anxiety and will admit to hating mathematics, when we give them space to explore mathematics and bring their existing knowledge to the problem, they can make some amazing mathematical discoveries.  Inquiry-based learning (IBL) is perfect structure to provide these type of opportunities. In this paper, we will examine one inquiry-based investigation in which MLA students were given the space to look for patterns which resulted in some original discoveries.
\end{abstract}

\section{INTRODUCTION}\label{Sect:Intro}

Mathematics for Liberal Arts (MLA) students will often claim that they hate mathematics and that they are not good at it.    Inquiry-based learning (IBL) provides a perfect opportunity for students to change that narrative as it allows them to use creativity and their own strengths to solve problems.  I encountered such a situation in my MLA classes recently when my students were exploring the \textbf{Rascal Triangle} \cite{ALT}.   However, in order for my students to be in a position to be that successful, a classroom environment that supported inquiry needed to be established; this section describes the course and several factors that helped set the stage for my students discoveries.

\subsection{Mathematical Explorations}\label{Ss:course}  

At Westfield State, \emph{Mathematical Explorations}, our most popular MLA course, is populated by students who are typically from the humanities or social sciences and have very negative attitudes toward mathematics \cite{PH1}. One of the main goals in this course courses is to have our students experience mathematics as an artistic, creative and humanistic endeavor.  As a result, the topics in the course are primarily from pure mathematics, (although the specific choices are left up to the instructor), and the course is taught using a cooperative, inquiry-based learning approach. See \cite{VE} and \cite{JF2} for more details about our classroom and \cite{EFHV} for a series of freely-available inquiry-based books we have written for this type of course. 
\smallskip

The room in which most of our Mathematical Explorations classes are taught are equipped with nine tables at which the students work cooperatively in groups of 3 to 5, while the instructor moves about the room answering and asking questions, offering encouragement and occasionally  providing suggestions.  
\smallskip

During class time, the students are engaged in mathematical tasks, either working through a set of carefully constructed questions designed to guide the students to (re)discover the solution, or working on an open-ended big question.  It was this latter approach that led to my students' discoveries.

\subsection{Course Aspects}
\begin{description}
	\item[Sharing] Throughout the semester, students were asked to share ideas and patterns with the class.  This started with the $\boldsymbol{3a+5b}$ problem: What are all the possible values for $3a+5b$ when $a$ and $b$ are non-negative integers? See \cite{JF1} for more details about solutions to this problem.  The solutions to this problem were not original, but students were asked throughout the investigation to share patterns they had found and why they believed the patterns were true. This helped establish a culture of sharing ideas and the reasoning behind them, whether right or wrong, so by the time we got to the Rascal Triangle exploration, the sharing of ideas was expected.  As one student wrote, "Yet once we find a solution, something that actually works for the problem, it is exhilarating. We become so excited to explain this to the teacher and present it to the class."

	\item[\textit{The Proof}] We watched the PBS special, "The Proof" \cite{Pf} in class and the students were asked to read a portion of the Introduction to \textit{Modular elliptic curves and Fermat's Last Theorem} by Andrew Wiles \cite{AW} and identify verbs that corresponded to actions taken by Wiles while working on the Taniyama-Shimura conjecture.  The identification of these actions helped the students see that while the problems they work on are much simpler than Taniyama-Shimura, their process was similar to that undertaken by Wiles.  They also saw that mistakes could be valuable and not something to be feared.
	
	\item[Mixing the groups] I would change the groups regularly so that, as as one student put it, "\dots as classes went on, the entire class became friends, which made the experience [in this course] so much more enjoyable."
	Another student wrote:
	\begin{quote}
		One other thing that was great about the environment of the class was that I actually got to know my classmates. In some other strictly lecture based classes, I would often go a whole semester without even getting to know probably around 80\% of the class because there were no opportunities to do so. In the group work environment however, it allowed us to bond with each other because we were all working towards a common goal, and I definitely was friendly with a strong majority of the students in class by the end of the semester, and also made a few close friends.
	\end{quote}
	
	\item[Clapping Rhythms and Pascal's Triangle] Prior to exploring the Rascal Triangle, I had my students do a clapping and rhythm activity from \textit{Discovering the Art of Mathematics: Music} \cite[Chapter 2] {CvR} that resulted in an exploration of patterns in Pascal's Triangle.  One student wrote,
		\begin{quote}
			For instance, we found that when choosing a selected amount of numbers diagonally, the sum of those numbers will allow the creation of a "hockey stick" shape, as the sum is also a number within Pascal's triangle, making the formation possible. This activity was not about long, complicated algorithmic equations, it was an act of discovery, observation and can even be seen as art.
		\end{quote}
\end{description}

These factors, and, of course, the students themselves, created an environment that led to the discoveries described in \secref{S:studentpatterns}.

\section{THE RASCAL TRIANGLE}\label{S:RasTri}
This section contains a brief review of the Rascal Triangle if the reader is unfamiliar with it.

\begin{figure}[H]
		\centering
		\includegraphics[scale=0.5]{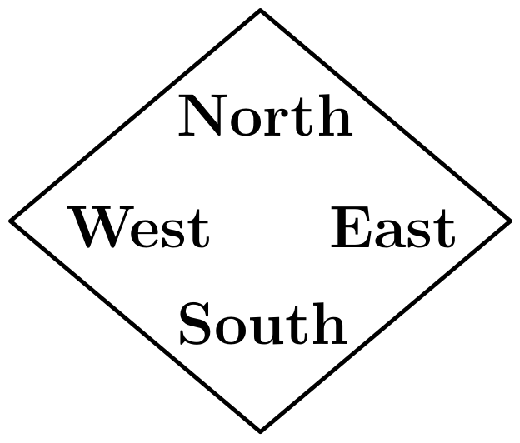}
		\caption{\textbf{North}, \textbf{South}, \textbf{East} and \textbf{West} entries in a triangular array.}\label{Fi:nsew}
	\end{figure}

 In 2010, middle school students Alif Anggaro, Eddy Liu and Angus Tulloch were asked to determine the next row of numbers in the following triangular array:
	\begin{figure}[H]
		\centering
		\includegraphics[scale=0.5]{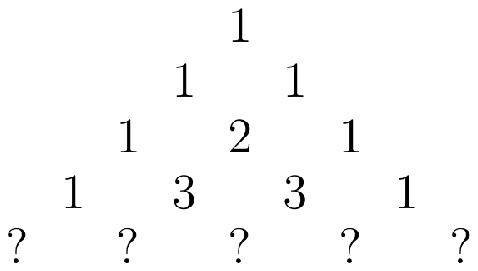}
		\caption{A triangular array.}\label{Fi:RasTriStart}
	\end{figure}
Instead of providing the row from Pascal's Triangle that the instructor expected,
	\begin{center}
		1\qquad 4\qquad 6\qquad 4\qquad 1
	\end{center}
they produced the row
	\begin{center}
		1\qquad 4\qquad 5\qquad 4\qquad 1.
	\end{center}
They did this by using the rule that the outside numbers are 1s and the inside numbers are determined by the \textbf{diamond formula}
	\begin{equation*}\label{Eq:ALT}
		\mathbf{ South} = \dfrac{\mathbf{East}\cdot \mathbf{West} + 1}{\mathbf{North}}
	\end{equation*}
where \textbf{North}, \textbf{South}, \textbf{East} and \textbf{West} form a diamond in the triangular array as in \figref{Fi:nsew}.
	
Continuing with this rule Anggaro, Liu and Tulloch created a number triangle they called the \textbf{Rascal Triangle} \cite{ALT}.
	\begin{figure}[H]
		\centering
		\includegraphics[scale=0.5]{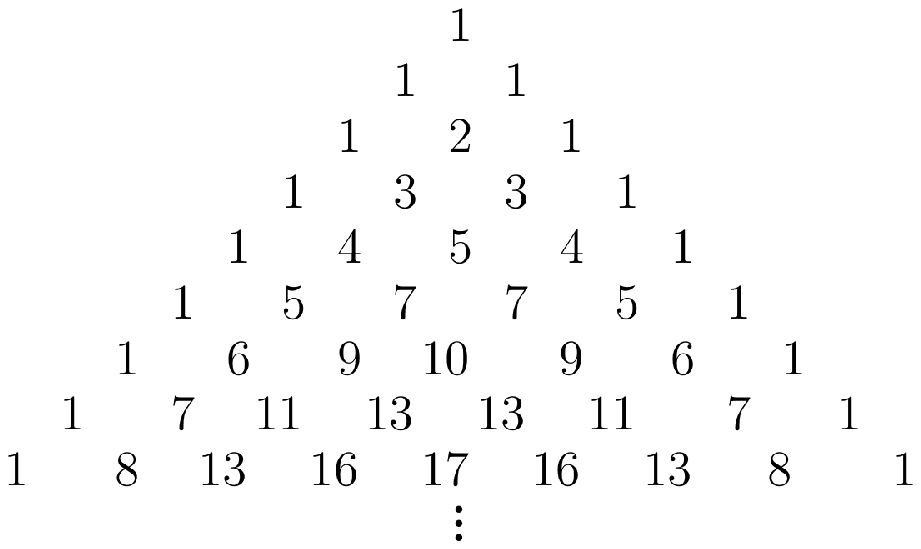}
		\caption{The Rascal Triangle.}\label{Fi:RasTri}
	\end{figure}

Because the diamond formula involves division, their instructor challenged Anggaro, Liu and Tulloch to prove that it would always result in an integer.  They did this by using the fact that the diagonals in the Rascal Triangle formed arithmetic sequences; see \cite{ALT} for details.  

In the Spring 2015 semester, students in a \emph{Mathematical Explorations} classes taught by one of my colleagues, Julian Fleron, discovered that the Rascal Triangle can also be generated by the rule that the outside numbers are 1s and the inside numbers are determined by the formula
	\begin{equation*}\label{Eq:JFclass}
		\mathbf{South} = \mathbf{East}+ \mathbf{West} -\mathbf{North} + 1.
	\end{equation*}
	
This formula also follows from the arithmetic progressions along the diagonals. See \cite{JF3} for details about this discovery; which, as best as can be determined, was unknown at the time.  Thus, the Rascal Triangle has the property that for any diamond containing 4 entries, the \textbf{South} entry satisfies two equations: 
	\begin{align}
		\mathbf{South} &= \dfrac{\mathbf{East}\cdot \mathbf{West} + 1}{\mathbf{North}}\label{eq:RasForm1}\\
		\mathbf{South} &= \mathbf{East}+ \mathbf{West} -\mathbf{North} + 1\label{eq:RasForm2}
	\end{align}
	
The fact that both \twoeqref{eq:RasForm1}{eq:RasForm2} can be used to generate the Rascal Triangle was intriguing; and, as part of an effort to better understand the Rascal Triangle, I had one of my \emph{Mathematical Explorations} classes explore number triangles in general and the Rascal Triangle in particular.  What followed was remarkable; my students became incredibly engaged with the exploration and some of them produced some original results, which are described next.

 \section{STUDENT PATTERNS AND NUMBER TRIANGLES}\label{S:studentpatterns}

After the exploration of Pascal's Triangle, I gave my MLA class the first six rows of the Rascal Triangle and asked them to find a method for extending the triangle in a manner that was consistent with the first six rows.  The hope was that they would identify the arithmetic sequences on the diagonals and eventually (re)-discover  \twoeqandorref{eq:RasForm1}{eq:RasForm2}.  The only instruction given was the reminder that Pascal's Triangle could be generated by the rule the \textbf{South} = \textbf{East} + \textbf{West}, so students were free to find any rule that made sense to them.  Students found this to be liberating; as one student wrote in her final refection (as a letter to a future student):
\begin{quote}
\dots Professor Hotchkiss gave us certain things, such as a number triangle, to find patterns in, every table's group would dive in and get so into what they were trying to figure out. We ended up figuring out crazy patterns, formulas, and equations in all kinds of different mathematical things that were presented to us in class. The best part was, Professor Hotchkiss never told us there was ever going to be a very certain correct answer. He told us to work through it and discover. Discover\dots That's a word I never thought I would here [sic] in math. I always thought of math as a very set way of doing things with only right or wrong answers. I always thought there was no discovering new things in math, you just learned what was already established in the text book. 

\dots We, as a class, literally found new mathematic patterns in a new number triangle that have never been discovered before. How cool is that? Pretty cool, I know. 
\end{quote}

Another student wrote
\begin{quote}
	\dots instead of being told how to  "do" the number triangles, Hotchkiss encouraged us to look for patterns in them. Not a specific pattern, ANY pattern. This was vital because it put the students in a position of power and gave them a thirst for finding more patterns. Once we got frustrated with finding something new or hopeless with a specific "mini-prompt" we would get a bit of assistance from Hotchkiss that would sort of "nudge us along" if you will, almost like a dad giving his son a little push while learning how to ride a bike. 
\end{quote}

\subsection{Student Patterns in the Rascal Triangle}\label{Ss:StudRasTriPats} Following "discovery" of the Rascal Triangle via the arithmetic sequences on the diagonals, I asked my students to find other rules or patterns that would generate the Rascal Triangle.  \twoeqref{eq:RasForm1}{eq:RasForm2} had not yet been introduced, and as before, very little instructions were given so students had the freedom to derive any patterns that made sense to them.

\begin{figure}[H]
		\centering
		\includegraphics[scale=0.5]{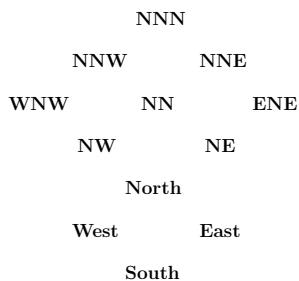}
		\caption{Location of Entries}\label{Fi:DirMap}
	\end{figure}

The patterns that were presented to the class were described using an informal notation involving directions relative to the the South entry, as shown in \figref{Fi:DirMap}, that was based on Angarro, Liu and Tulloch's diamond description of the entries in the Rascal Triangle in \cite{ALT}. 

The following definition was also helpful.
\begin{Def}
	For any number triangle, the diagonals running from right to left are called the \textbf{major} diagonals while the diagonals running from left to right are called the \textbf{minor} diagonals.
\end{Def}

	\begin{figure}[H]
		\centering
		\begin{subfigure}[b]{0.2\textwidth}
			\includegraphics[width=0.675\textwidth]{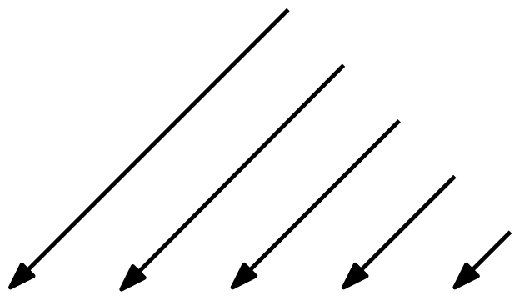}
			\subcaption{major diagonals.}\label{Fi:majord}
		\end{subfigure}
		\hskip0.5in
		\begin{subfigure}[b]{0.2\textwidth}
			\includegraphics[width=0.675\textwidth]{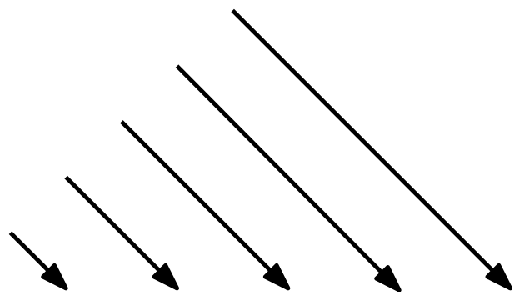}
			\subcaption{Minor diagonals.}\label{Fi:minord}
		\end{subfigure}
	\end{figure}

Due to time constraints, proofs of the student patterns were not presented in class. However, they are accessible to undergraduate mathematics majors, and I would expect that they would be accessible to most MLA students as well. Proofs of these patterns (in a more generalized setting) can be found in \cite{PH3}.  

Several groups provided methods for extending the number triangle, but two of these were of particular interest.  The first pattern was called the \textbf{T-Meg Rule}.
\begin{enumerate}
	\setcounter{enumi}{\value{temp}}
	\item\label{T-Meg} \underbar{T-Meg Rule}: Tim and Meaghan observed that starting with row 4 of the Rascal Triangle, we could determine \textbf{South} from \textbf{North} and the first two entries, $\mathbf{North}_0$ and $\mathbf{North}_1$, on \textbf{North}'s row.  That is,
		\begin{equation*}
			\mathbf{South} =  \mathbf{North} + \mathbf{North}_0 + \mathbf{North}_1.
		\end{equation*}
		\begin{figure}[H]
			\centering
			\includegraphics[scale=0.5]{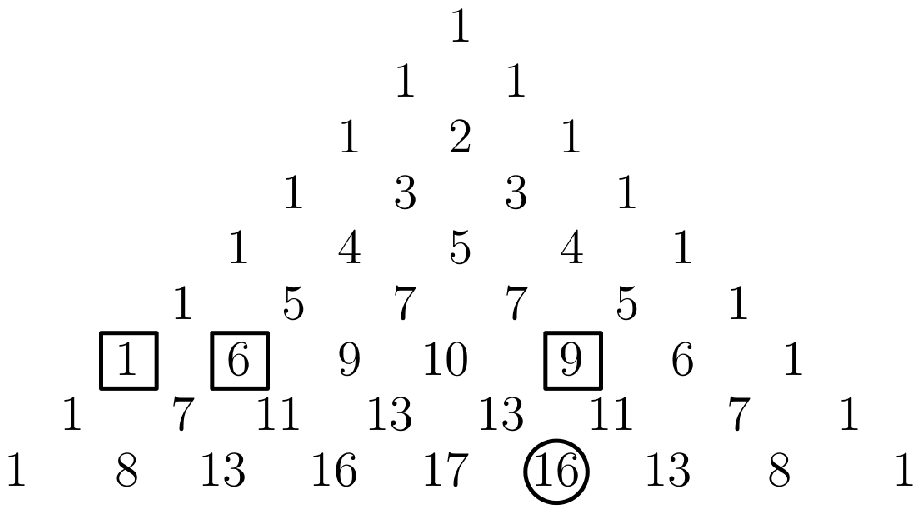}
			\caption{T-Meg Rule.}
		\end{figure}
	For example,
		\begin{equation*}
			16 = \mathbf{North} + \mathbf{North}_0 + \mathbf{North}_1 = 9 +1+6
		\end{equation*}
		
	Note that T-Meg's Rule also holds for the second entry in a row as well.
		\begin{equation*}
			8 = \mathbf{North} + \mathbf{North}_0 + \mathbf{North}_1 = 1+1+6
		\end{equation*}
	\setcounter{temp}{\value{enumi}}
\end{enumerate}
As Meaghan wrote about her role in finding the T-Meg Rule (\exref{T-Meg} below)
	\begin{quote}
		The hockey still pattern in Pascal's triangle is what helped me find the T-Meg pattern in Rascal's triangle. I knew there was a similar hockey stick pattern in Pascal's triangle which made me look for other ways addition was used in Rascals Triangle. This is how I was able to find the T-Meg pattern. That was a huge turning point for me because I had found a pattern in the triangle that no one else had previously found. That definitely made me more confident in my mathematical abilities.
	\end{quote}
	
\begin{enumerate}
	\setcounter{enumi}{\value{temp}}
	\item\label{Ashley} \underbar{Ashley's Rule}: Ashley observed that we could determine \textbf{South} from \textbf{West}, \textbf{East},  \textbf{NW} plus a "diagonal factor". The factor was associated with the minor diagonal containing the \textbf{South} entry, and was equal to $2-k$, where $k=0$ corresponds to the minor diagonal consisting of all 1s.
  		\begin{figure}[H]
			\centering
			\includegraphics[scale=0.5]{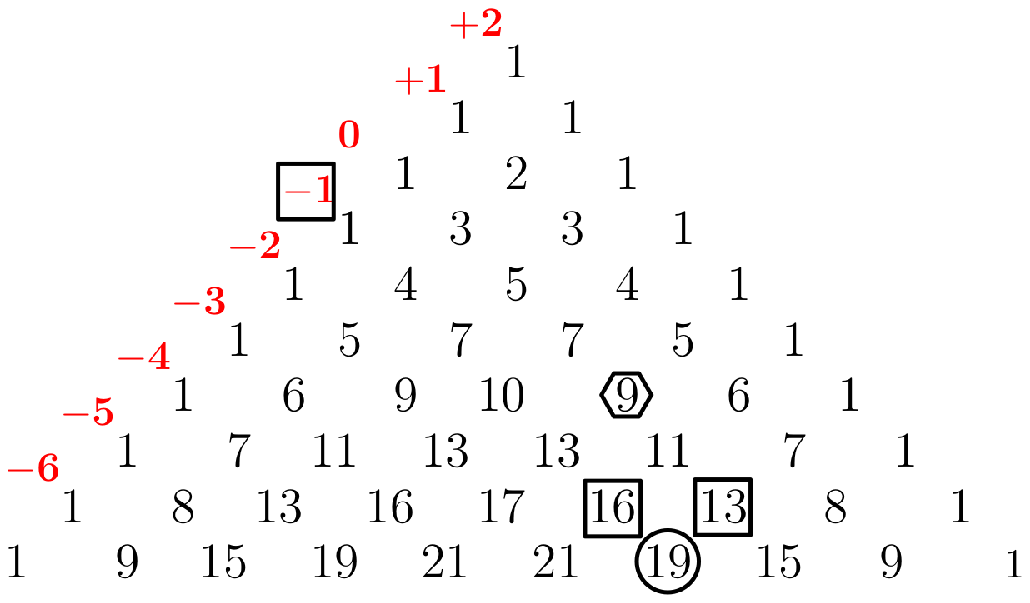}
			\caption{Ashley's Rule}
		\end{figure}
	For example, 
		\begin{equation*}
			19 = \mathbf{West} + \mathbf{East} -\mathbf{NW} - \mathbf{diagonal\ factor} = 16+13-9-1
		\end{equation*}
	\setcounter{temp}{\value{enumi}}
\end{enumerate}

While looking for methods of generating the Rascal Triangle, one student observed two patterns involving diamonds within the Rascal Triangle.

\begin{enumerate}
	\setcounter{enumi}{\value{temp}}	
	\item\label{John} \underbar{John's Diamond Patterns}: 
		\begin{enumerate}
			\item \underbar{John's Odd Diamond Pattern}: The average of the $8n$ numbers that form the edge of a diamond in the Rascal Triangle with $2n+1$ numbers on each side is the number in the middle of the diamond.  For example, for the diamond in \figref{Fi:OddDiamond},
				\begin{figure}[H]
					\centering
					\includegraphics[scale=0.5]{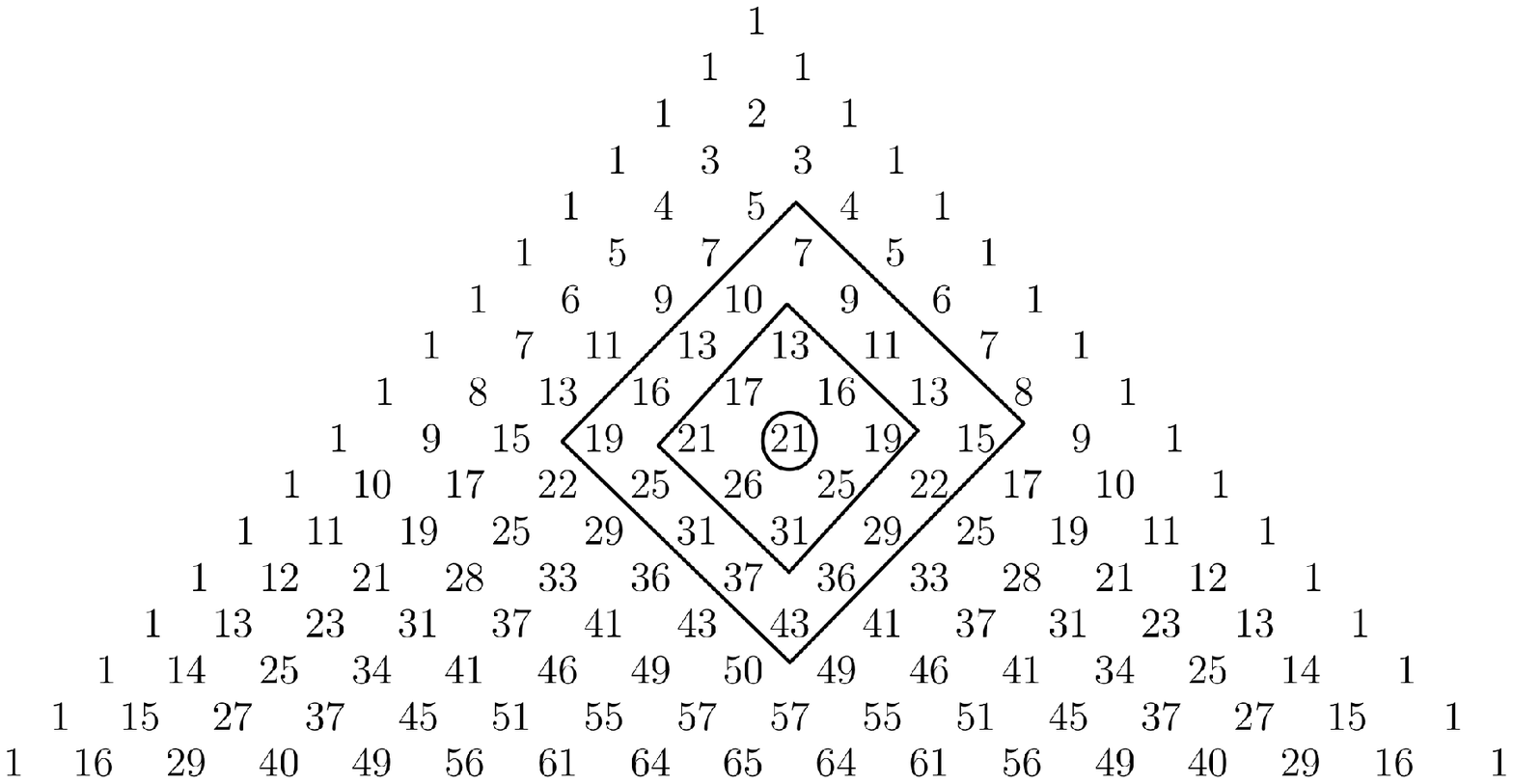}
					\caption { John's Odd Diamond Pattern.}\label{Fi:OddDiamond}
				\end{figure}

				\begin{align*}
					&\dfrac{13+16+19+25+31+26+21+17}8 = \dfrac{168}8 = 21\\
					&{\textstyle\frac{7+9+11+13+15+22+29+36+43+37+31+25+19+16+13+10}{16} = \frac{336}{16} = 21}
				\end{align*}
	
			\item \underbar{John's Even Diamond Pattern}: If you form a $2\times2$-diamond in the Rascal Triangle and a diamond with $2n$ numbers on each side that has the $2\times2$-diamond in the center, then the average of the $8n-4$ numbers along the edges of the outer diamond is equal to the average of the 4 numbers along the edge of the $2\times2$-diamond. For example, in \figref{Fi:EvenDiamond} 
				\begin{figure}[H]
					\centering
					\includegraphics[scale=0.5]{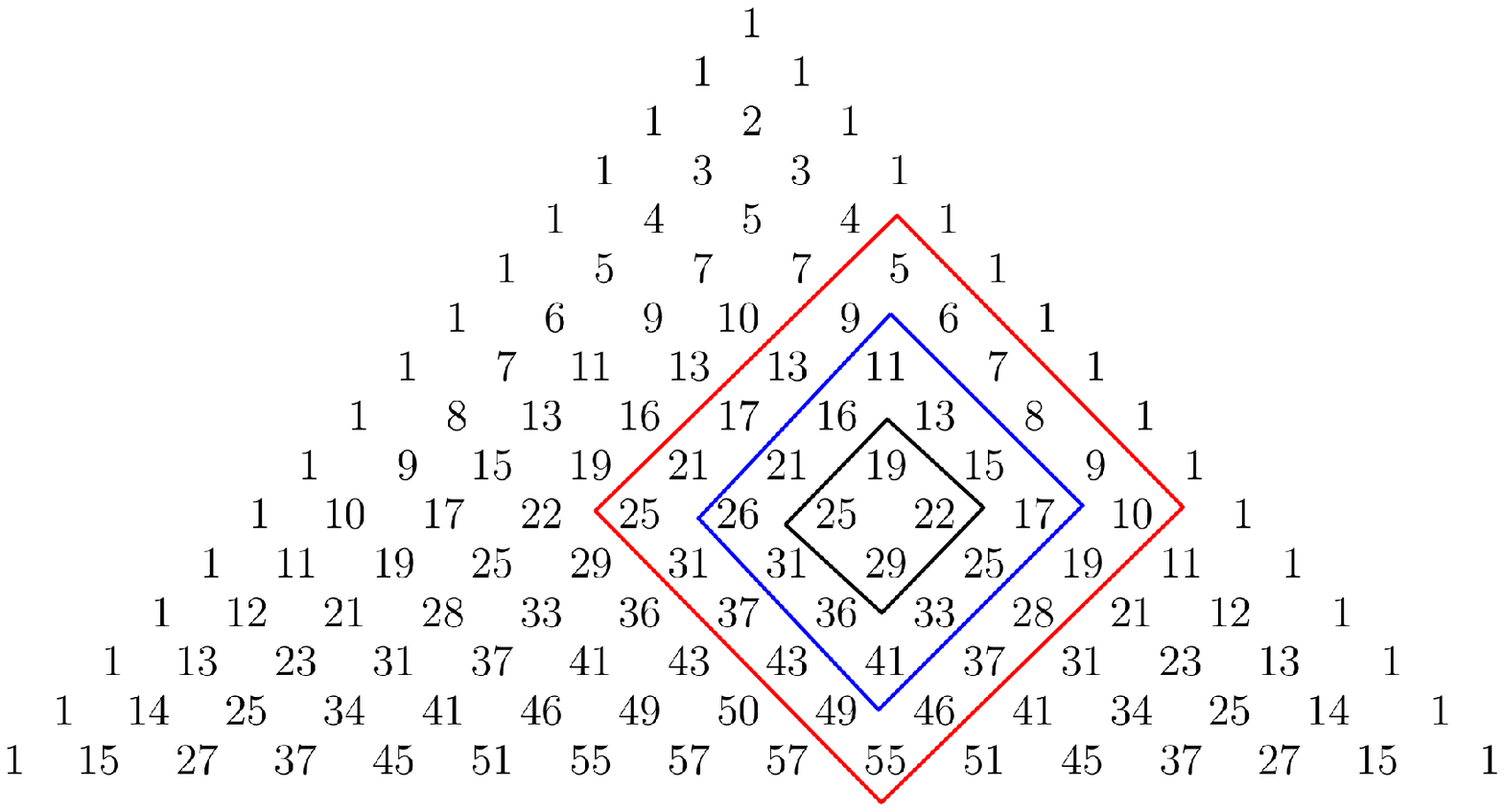}
					\caption{John's Even Diamond pattern}\label{Fi:EvenDiamond}
				\end{figure}

				\begin{align*}
					&\dfrac{19+22+29+25}4\\ 
					&= \dfrac{95}4 =23.75;\\
					&\dfrac{11+13+15+17+25+33+41+36+31+26+21+16}{12}\\ 
					&= \dfrac{285}{12} = 23.75;\\
					&{\textstyle\frac{5+6+7+8+9+10+19+28+37+46+55+49+43+37+31+25+21+17+13+9}{20}}\\ 
					&= \dfrac{475}{20}=23.75.
					\intertext{John's original formulation of this pattern was that the the sum of the $8n-4$ numbers along the edge of the outer diamond was $2n-1$ times the sum of the 4 numbers along the inner diamond.}
					&19+22+29+25\\ 
					&= 95;\\
					&11+13+15+17+25+33+41+36+31+26+21+16\\ 
					&= 285 = 3\cdot 95;\\
					&{\scriptstyle 5+6+7+8+9+10+19+28+37+46+55+49+43+37+31+25+21+17+13+9}\\ 
					&= 475 = 5\cdot 95.
				\end{align*}
		\end{enumerate}
\setcounter{temp}{\value{enumi}}
\end{enumerate}

\section{CONCLUSION}

The results in this paper are the results of explorations by Mathematics for Liberal Arts students looking for patterns in the Rascal Triangle.  In particular, the discovery that the Rascal Triangle could be generated by a several different addition rules, illustrate the power of "non-expert" eyes seeing patterns that the experts could not see because they already knew the patterns.   One student wrote at the end of the semester
\begin{quote}
One moment when I was inspired by math in this class was when we were discussing the Rascalls [sic] triangle. It amazed me how the students in my class whom we all worked together were able to find proofs on things that very high level mathematicians had been looking at for years. I had never once been in an environment in a Math classroom in which I felt like everyone there was willing to help and support each other. We were not treated like cogs in a machine that we had no control of. We were treated like actually able bodied human beings and were able to look at things in ways that were not controlled by the teacher. 
\end{quote}
\medskip

Theses student investigations inspired me to examine the algebraic structure of the Rascal Triangle more deeply and resulted in the exploration of \textit{Generalized Rascal Triangles}, see \cite{PH3} for details.

\section*{ACKNOWLEGEMENTS}
The author would like to acknowledge his students who observed the patterns discussed in this paper: Ashley Craig, John Coulombe, Timothy Schreiner, Meaghan Sparks and Evan Wilson.


\begin{thebibliography}{10}
	\bibitem{ALT}
	Anggoro, A., Liu, E., and Tulloch, A.,
	The Rascal Triangle,
	\emph{The College Mathematics Journal},
	\textbf{41}, No. 5, November 2010, pp. 393-395.
	
	\bibitem{CBMS}
	Blair, R., Kirkman, E.E., Maxwell, J.W., 
	Statistical abstract undergraduate programs in the mathematical sciences in the United States: 2015 CBMS survey,
	American Mathematical Society, Providence, RI, 2018. 
	
	\bibitem{VE}
	Ecke, V.,
	Our Inquiry-Based Classroom,
	retrieved from \url{http://www.artofmathematics.org/blogs/vecke/our-inquiry-based-classroom}.
	
	\bibitem{EFHV} 
	Ecke, V., Fleron, J., Hotchkiss, P. and von Renesse, C, 
	Books: Inquiry-based Learning Guides, 
	retrieved from \url{http://www.artofmathematics.org/books}.
	
	\bibitem{JF1}
	Fleron, J.,
	$3a+5b$ Proofs,
	retrieved from \url{http://www.artofmathematics.org/blogs/jfleron/3a5b-proofs}.
	
	\bibitem{JF2}
	Fleron, J.,
	Classroom Vignette
	retrieved from \url{http://www.artofmathematics.org/blogs/jfleron/classroom-vignette}.
	
	\bibitem{JF3}
	Fleron, J., 
	Fresh Perspectives Bring Discoveries,
	\emph{Math Horizons}, \textbf{24}, No. 3, February 2017, p. 15.
	
	\bibitem{OEIS}
	Online Encyclopedia of Integer Sequences
	retrieved from \url{https://oeis.org/}
	
	\bibitem{PH1}
	Hotchkiss, P., 
	Audience: Learning about our MLA Students, 
	retrieved from \url{http://www.artofmathematics.org/blogs/photchkiss/audience-learning-about-our-mla-students}.
	
	\bibitem{PH2}
	Hotchkiss, P.,
	Movie "Proof": How our students view the process of mathematics,
	retrieved from \url{http://www.artofmathematics.org/blogs/photchkiss/movie-the-proof-how-our-students-view-the-process-of-mathematics}.
	
	\bibitem{PH3}
	Hotchkiss, P
	Generalized Rascal Triangles,
	preprint, 2019.
	
	\bibitem{Pf}
	"The Proof", \emph{NOVA}, WGBH, Boston, 1997. DVD 
	
	\bibitem{CvR}
	von Renesse, C. with Fleron, J. Hotchkiss, P. and Ecke, V., \emph{Discovering the Art of Mathematics: Music}, \url{http://www.artofmathematics.org/books/music}, 2015.
	
	\bibitem{AW}
	Wiles, A.,
	Modular, elliptic curves and Fermat's Last Theorem
	\emph{Annals of Mathematics}, \textbf{141}, No. 3, 1995, pp. 443-551
\end{thebibliography}
\end{document}